\documentclass{article}
\usepackage{amsfonts}
\usepackage{amsmath}

\input{tcilatex}
\begin{document}

\begin{titlepage}
\title{\bf Conformal Euler-Lagrangian Equations on Para-Quaternionic K\"{a}hler Manifolds}
\author{Zeki Kasap {\footnote{Corresponding Author. E-mail: zkasap@pau.edu.tr; 
Tel:+9 0 258 296 11 82; Fax:+9 0 258 296 12 00}} \\
  {\small Department of Elementary Education, Faculty of Education,  Pamukkale University}
\\{\small  20070 Denizli, Turkey.}
\\ Mehmet Tekkoyun \footnote{Author. E-mail: tekkoyun@pau.edu.tr; Tel: +9 0 258 296 36 16} 
\\ {\small Department of Mathematics, Pamukkale University}
\\{\small  20070 Denizli, Turkey.} }
\date{\today}
\maketitle
\begin{abstract}
This manuscript presents an attempt to introduce Lagrangian formalism for mechanical systems using
para-quaternionic K\"{a}hler  manifolds, which represent an interesting multidisciplinary field of research. 
In addition to, the
geometrical-physical results related to para-quaternionic K\"{a}hler mechanical systems are also given.



\end{abstract}
\end{titlepage}

\section{Introduction}

Today, many branches of science are into our lives. One the branches is
mathematics that has multiple applications. In particular, differential
geometry and mathematical physics have a lots of different applications.
These applications are used in many areas. One of them are on geodesics.
Geodesics are known the shortest route between two points. Time-dependent
equations of geodesics can be easily found with the help of the
Euler-Lagrange equations. We\ can say that differential geometry provides a
good working area for studying Lagrangians of classical mechanics\ and field
theory. The dynamic equation for moving bodies is obtained for Lagrangian
mechanic. These dynamic equation is illustrated as follows:

\textbf{Lagrange Dynamics Equation} \cite{klein,deleon,abraham}:\textbf{\ }%
Let $M$ be an $n$-dimensional manifold and $TM$ its tangent bundle with
canonical projection $\tau _{M}:TM\rightarrow M$. $TM$ is called the phase
space of velocities of the base manifold $M$. Let $L:TM\rightarrow R$ be a
differentiable function on $TM$\ called the \textbf{Lagrangian function}. We
consider the closed 2-form on $TM$ given by $\Phi _{L}=-dd_{J}L$ ( if $%
J^{2}=-I,$ $J$\ is a complex structure and if $J^{2}=I$, $J$ is a
paracomplex structure, $Tr(J)=0$ ). Consider the equation 
\begin{equation}
i_{X}\Phi _{L}=dE_{L}.  \label{1.1}
\end{equation}%
Then $X$ is a vector field, we shall see that (\ref{1.1}) under a certain
condition on $X$ is the intrinsical expression of the Euler-Lagrange
equations of motion. This equation is named \ as \textbf{Lagrange dynamical
equation}. We shall see that for motion in a potential, $E_{L}=V\left(
L\right) -L$\ \ is an energy function and $V=J(X)$\ a Liouville vector
field. Here $dE_{L}$ denotes the differential of $E$. The triple $(TM,\Phi
_{L},X)$ is known as \textbf{Lagrangian system} on the tangent bundle $TM$.
If it is continued the operations on (\ref{1.1}) for\ any coordinate system $%
(q^{i}(t),p_{i}(t))$, infinite dimension \textbf{Lagrange's\ equation is }%
obtained the form below: 
\begin{equation}
\frac{dq^{i}}{dt}=\dot{q}^{i}\ ,\ \frac{d}{dt}\left( \frac{\partial L}{%
\partial \dot{q}^{i}}\right) =\frac{\partial L}{\partial q^{i}}\ ,\
i=1,...,n.  \label{1.2}
\end{equation}%
There are many studies about Lagrangian dynamics, mechanics, formalisms,
systems and equations (see detail \cite{deleon1}). There are real, complex,
paracomplex and other analogues. It is well-known that Lagrangian analogues
are very important tools. They have a simple method to describe the model
for mechanical systems. The models about mechanical systems are given as
follows.

Some examples of the Lagrangian is applied to model the problems include
harmonic oscillator, charge $Q$ in electromagnetic fields, Kepler problem of
the earth in orbit around the sun, pendulum, molecular and fluid dynamics, $%
LC$ networks, Atwood's machine, symmetric top etc. Let's remember some work
done. \textit{Vries} shown that the Lagrangian motion equations have a very
simple interpretation in relativistic quantum mechanics \cite{vries}.
Paracomplex analogue of the Euler-Lagrange equations was obtained in the
framework of para-K\"{a}hlerian manifold and the geometric results on a
paracomplex mechanical systems were found by \textit{Tekkoyun} \cite%
{tekkoyun2005}. Electronic origins, molecular dynamics simulations,
computational nanomechanics, multiscale modelling of materials fields were
contributed by \textit{Liu} \cite{liu}. Bi-paracomplex analogue of
Lagrangian systems was shown on Lagrangian distributions by \textit{Tekkoyun 
}and\textit{\ Sari} \cite{tekkoyun2010}. \textit{Tekkoyun }and\textit{\ Yayli%
} presented generalized-quaternionic K\"{a}hlerian analogue of Lagrangian
and Hamiltonian mechanical systems. Eventually, the geometric-physical
results related to generalized-quaternionic K\"{a}hlerian mechanical systems
are provided \cite{tekkoyun2011}.

Nowadays, there are many studies about Euler-Lagrangian dynamics, mechanics,
formalisms, systems and equations \cite{deleon, deleon1, maxim, tekkoyun,
tekkoyun1} and there in. There are real, complex, paracomplex and other
analogues. As known it is possible to produce different analogous in
different spaces. Quaternions were invented by Sir William Rowan \textbf{%
Hamiltonian }as an extension to the complex numbers. \textbf{Hamiltonian's }%
defining relation is most succinctly written as:%
\begin{equation*}
i^{2}=j^{2}=k^{2}=-1,\text{~}ijk=-1.
\end{equation*}%
Split quaternions are given by%
\begin{equation*}
i^{2}=-1,\text{ }j^{2}=1=k^{2},\text{ }ijk=1.
\end{equation*}%
Generalized quaternions are defined as%
\begin{equation*}
i^{2}=-a,\text{ }j^{2}=-b,\text{ }k^{2}=-ab,\text{ }ijk=-ab.
\end{equation*}%
If it is compared to the calculus of vectors, quaternions have slipped into
the realm of obscurity. They do however still find use in the computation of
rotations. Lots of physical laws in classical, relativistic, and quantum
mechanics can be written pleasantly by means of quaternions. Some physicists
hope they will find deeper understanding of the universe by restating basic
principles in terms of quaternion algebra \cite{dan, dancer, kula1, kula2,
ata1, ata2, jafari}.

In the present paper, we present equations related to Lagrangian mechanical
systems on a generalized-quaternionic K\"{a}hler manifold.

\section{Preliminaries}

In this study, all the manifolds and geometric objects are $C^{\infty }$\
and the Einstein summation convention is in use. Also, $A$, $F(TM)$, $\chi
(TM)$\ and $\Lambda ^{1}(TM)$\ denote the set of paracomplex numbers, the
set of (para)-complex functions on $TM$, the set of (para)-complex vector
fields on $TM$\ and the set of (para)-complex 1-forms on $TM$, respectively.
The definitions and geometric structures on the differential manifold $M$\
given in \cite{cruceanu} may be extended to $TM$\ as follows:

\section{Conformal Geometry}

In mathematics, a conformal map is a function which preserves angles. In the
most common case the function is between domains in the complex plane.
Conformal maps can be defined between domains in higher dimensional
Euclidean spaces, and more generally on a Riemann or semi-Riemann manifold.
Conformal geometry is the study of the set of angle-preserving (conformal)
transformations on a space. In two real dimensions, conformal geometry is
precisely the geometry of Riemann surfaces. In more than two dimensions,
conformal geometry may refer either to the study of conformal
transformations of "flat" spaces (such as Euclidean spaces or spheres), or,
more commonly, to the study of conformal manifolds which are Riemann or
pseudo-Riemann manifolds with a class of metrics defined up to scale. A
conformal manifold is a differentiable manifold equipped with an equivalence
class of (pseudo) Riemann metric tensors, in which two metrics $g^{\prime }$
and $g$ are equivalent if and only if%
\begin{equation}
g^{\prime }=\lambda ^{2}g  \label{1.3}
\end{equation}%
where $\lambda >0$ is a smooth positive function. An equivalence class of
such metrics is known as a conformal metric or conformal class \cite{wiki}.

\section{\textbf{Conformal} S\textbf{tructure}}

The linear distance between two points can be found easily by Reimann
metric, which is very useful and is defined inner product. Many scientists
have used the Riemann metric. Einstein was one of the first studies in this
field. Einstein discovered which the Riemannian geometry and successfully
used \ it to describe General Relativity in the 1910 that is actually a
classical theory for gravitation. However, the universe is really completely
not like Riemannian geometry. Each path between two points is not always
linear. Also, orbits of move objects may change during movement. So, each
two points in space may not be linear geodesic and need not\ to be.
Therefore, new metric\ is needed for non-linear distances like spherical
surface. Then, a method is required for converting nonlinear distance to
linear distance. Weyl introduced a metric with a conformal transformation in
1918.

\textbf{Definition 1.}\textit{\ Let }$M$\textit{\ an }$n$\textit{%
-dimensional smooth manifold.\ A\textbf{\ conformal structure} on }$M$%
\textit{\ is an equivalence class }$G$\textit{\ of Riemann metrics on }$M$%
\textit{. }A manifold with a conformal structure is called a \textbf{%
conformal manifold}.

\textbf{(i) }Two Riemann metrics $g$ and $g^{\prime }$ on $M$ \ are said to
be equivalent if and only if%
\begin{equation}
g^{\prime }=e^{\lambda }g  \label{1.4}
\end{equation}%
where $\lambda $ is a smooth function on $M$. The equation given by (\ref%
{1.4}) is called a \textbf{conformal structure}

\textbf{(ii)} A \textbf{Weyl structure} on $M$ is a map $F:G\rightarrow
\wedge ^{1}M$ satisfying%
\begin{equation}
F(e^{\lambda }g)=F(g)-d\lambda  \label{1.5}
\end{equation}%
where $G$ \ is a conformal structure. Note that a Riemann metric $g$ and a
one-form $\varphi $ determine a Weyl structure, namely $F:G\rightarrow
\wedge ^{1}M$\textit{\ where }$G$ is the equivalence class of $g$ and\textit{%
\ }$\ F(e^{\lambda }g)=\varphi -d\lambda $.

\textbf{Theorem 1.\ }\textit{A connection on the metric bundle }$\varphi $%
\textit{\ of a conformal manifold }$M$\textit{\ naturally induces a map }$%
F:G\rightarrow \wedge ^{1}M$\textit{\ } and (\ref{1.5}),\textit{\ and
conversely. Parallel translation of points in }$\varphi $\textit{\ by the
connection is the same as their translation by }$F$\textit{\ }\cite{folland}.

\section{ Generalized-Quaternionic K\"{a}hler Manifolds}

A generalized almost quaternion structure on the manifold $M$\ is a
subbundle of the bundle of endomorphisms of the tangent bundle of $M$, whose
standard fibre is the algebra of generalized quaternions. A generalized
almost quaternion structure on a pseudo-Riemannian manifold is called a
generalized quaternion-Hermitian if the following conditions hold:

i) The endomorphisms $F,G$\ and $H$\ of $T_{x}M$\ satisfy%
\begin{equation}
F^{2}=-aI,\text{ }G^{2}=-bI,~\text{ }H^{2}=-abI,\text{ }FG=H,GH=bF,HF=aG,
\label{2.1}
\end{equation}%
ii) The compatibility equations are given by, for $X,Y\in T_{x}M,$%
\begin{equation}
g(FX,FY)=ag(X,Y),\text{ }g(GX,GY)=bg(X,Y),\text{ }g(HX,HY)=abg(X,Y),
\label{2.2}
\end{equation}%
where $I$\ denotes the identity tensor of type (1,1) in $M.$\ In particular,
2-form $Q$\ defined by $Q(X,Y)=(X,FY)=(X,GY)=(X,HY)$\ on $M$\ is called the K%
\"{a}hler form of the endomorphisms $F,G$\ and $H$. If the K\"{a}hler form $%
Q $\ on $M$\ is closed, i.e. $dQ=0,$\ the manifold $M$\ is called a
generalized-quaternionic K\"{a}hler manifold \cite{arsen}. If $a=b=1$, $M$\
is quaternion manifold. If $a=1,b=-1$, $M$\ is para-quaternion manifold. The
bundle $V$\ is a set that locally admits a basis $\{F,G,H\}$\ satisfying (%
\ref{2.1}) and (\ref{2.2}) in any coordinate neighborhood $U\subset M$\ such
that $M=\cup U$\ \cite{dancer}. Then $V$ is called a
generalized-quaternionic structure in $M$. The pair $(M,V)$ denotes a
generalized-quaternionic manifold with $V$. The structure $V$ with such a
Riemannian metric $g$ is called a generalized-quaternionic metric structure$%
. $ The triple $(M,g,V)$ denotes a generalized-quaternion metric manifold.
Let $\left\{ x_{i},x_{n+i},x_{2n+i},x_{3n+i}\right\} ,$ $i=\overline{1,n}$
be a real coordinate system on a neighborhood $U$ of $M,$ and let $\left\{ 
\frac{\partial }{\partial x_{i}},\frac{\partial }{\partial x_{n+i}},\frac{%
\partial }{\partial x_{2n+i}},\frac{\partial }{\partial x_{3n+i}}\right\} $
and $\{dx_{i},dx_{n+i},dx_{2n+i},dx_{3n+i}\}$ be natural bases over $R$ of
the tangent space $T(M)$ and the cotangent space $T^{\ast }(M)$ of $M,$
respectively$.$ Taking into consideration (\ref{2.1}), then we can obtain
the expressions as follows:%
\begin{equation}
\begin{array}{c}
F(\frac{\partial }{\partial x_{i}})=a\frac{\partial }{\partial x_{n+i}},%
\text{ }F(\frac{\partial }{\partial x_{n+i}})=-a\frac{\partial }{\partial
x_{i}},\text{ }F(\frac{\partial }{\partial x_{2n+i}})=a\frac{\partial }{%
\partial x_{3n+i}},\text{ }F(\frac{\partial }{\partial x_{3n+i}})=-a\frac{%
\partial }{\partial x_{2n+i}} \\ 
G(\frac{\partial }{\partial x_{i}})=-b\frac{\partial }{\partial x_{2n+i}},%
\text{ }G(\frac{\partial }{\partial x_{n+i}})=b\frac{\partial }{\partial
x_{3n+i}},\text{ }G(\frac{\partial }{\partial x_{2n+i}})=-b\frac{\partial }{%
\partial x_{i}},\text{ }G(\frac{\partial }{\partial x_{3n+i}}\mathbf{)=}b%
\frac{\partial }{\partial x_{n+i}} \\ 
H(\frac{\partial }{\partial x_{i}})=-ab\frac{\partial }{\partial x_{3n+i}},%
\text{ }H(\frac{\partial }{\partial x_{n+i}})=-ab\frac{\partial }{\partial
x_{2n+i}},\text{ }H(\frac{\partial }{\partial x_{2n+i}})=-ab\frac{\partial }{%
\partial x_{n+i}},\text{ }H(\frac{\partial }{\partial x_{3n+i}})=-ab\frac{%
\partial }{\partial x_{i}}.%
\end{array}
\label{2.22}
\end{equation}

\section{ Generalized-Quaternionic Conformal K\"{a}hler Manifolds}

\textbf{Definition 2}.\textbf{\ }Let $(M,g,\nabla ,J_{\pm })$ be an almost
para/pseudo-Hermitian Weyl manifold. If $\ \nabla (J_{\pm })=0$, then one
says that this is a \textbf{(para)-K\"{a}hler Weyl manifold}. Note that
necessarily $J_{\pm }$ is integrable in this setting.

\textbf{Theorem 2 }.\textbf{\ }\textit{If }$(M,g,\nabla ,J_{\pm })$\textit{\
is a (para)-K\"{a}hler Weyl manifold with dimension }$n\geq 6$\textit{\ and
with }$H^{1}(M;R)=0$\textit{, then the underlying Weyl structure on }$M$%
\textit{\ is trivial.}

\textbf{Theorem 3.} \textit{If }$(M,g,\nabla ,J_{\pm })$\textit{\ is a
curvature (para)-K\"{a}hler Weyl manifold with dimension }$n\geq 6$\textit{\
and with }$H^{1}(M;R)=0$\textit{, then the underlying Weyl structure on }$M$%
\textit{\ is trivial.}

\textbf{Theorem} \textbf{4}. \textit{Let }$n\geq 6$\textit{. If }$(M,g,J\pm
,\nabla )$\textit{\ is a K\"{a}hler--Weyl structure, then the associated
Weyl structure is trivial, i.e. there is a conformally equivalent metric }$%
\tilde{g}=e^{2f}g$\textit{\ so that }$(M,\tilde{g},J\pm )$\textit{\ is K\"{a}%
hler and so that }$\nabla =\nabla ^{\tilde{g}}$ \cite%
{gilkey,gilkey2011,pedersen}.

After this part $W$ will be used instead of $J$. A manifold with a Weyl
structure is known as a Weyl manifold. $\lambda $ second structure was
chosen the minus sign. Because the condition of the structure required to
provide. $W_{\pm }^{2}=\pm Id$ \cite{miron}. If we rewrite (\ref{2.22})
equation with conformal structure, we obtain the following equations:%
\begin{equation}
\begin{array}{c}
W_{F}(\frac{\partial }{\partial x_{i}})=ae^{\lambda }\frac{\partial }{%
\partial x_{n+i}},\text{ }W_{F}(\frac{\partial }{\partial x_{n+i}}%
)=-ae^{-\lambda }\frac{\partial }{\partial x_{i}},\text{ }W_{F}(\frac{%
\partial }{\partial x_{2n+i}})=ae^{\lambda }\frac{\partial }{\partial
x_{3n+i}}, \\ 
W_{F}(\frac{\partial }{\partial x_{3n+i}})=-ae^{-\lambda }\frac{\partial }{%
\partial x_{2n+i}},\text{ }W_{G}(\frac{\partial }{\partial x_{i}}%
)=-be^{\lambda }\frac{\partial }{\partial x_{2n+i}},\text{ }W_{G}(\frac{%
\partial }{\partial x_{n+i}})=be^{\lambda }\frac{\partial }{\partial x_{3n+i}%
}, \\ 
W_{G}(\frac{\partial }{\partial x_{2n+i}})=-be^{-\lambda }\frac{\partial }{%
\partial x_{i}},\text{ }W_{G}(\frac{\partial }{\partial x_{3n+i}}%
)=be^{-\lambda }\frac{\partial }{\partial x_{n+i}},\text{,}W_{H}(\frac{%
\partial }{\partial x_{i}})=-abe^{\lambda }\frac{\partial }{\partial x_{3n+i}%
}, \\ 
W_{H}(\frac{\partial }{\partial x_{n+i}})=-abe^{\lambda }\frac{\partial }{%
\partial x_{2n+i}},\text{ }W_{H}(\frac{\partial }{\partial x_{2n+i}}%
)=-abe^{-\lambda }\frac{\partial }{\partial x_{n+i}},\text{ }W_{H}(\frac{%
\partial }{\partial x_{3n+i}})=-abe^{-\lambda }\frac{\partial }{\partial
x_{i}}\emph{.}%
\end{array}
\label{2.3}
\end{equation}%
We continue our studies thinking of the $(M,g,\nabla ,W_{\pm })$ instead of 
\textit{the almost para/pseudo-K\"{a}hler Weyl manifolds} $(M,g,\nabla
,J_{\pm })$.

\section{Conformal Euler-Lagrangian Mechanical Systems}

Here, we obtain Euler-Lagrange equations for quantum and classical mechanics
by means of a canonical local basis $\{F,G,H\}$ of $V$ that they defined on
a generalized-quaternionic K\"{a}hler manifold $(M,g,V).$ Firstly, let $F$
take a local basis element on the generalized-quaternionic K\"{a}hler
manifold $(M,g,V),$ and $\left\{ x_{i},x_{n+i},x_{2n+i},x_{3n+i}\right\} $
be its coordinate functions. Let semispray be the vector field $X$
determined by 
\begin{equation}
X=X^{i}\frac{\partial }{\partial x_{i}}+X^{n+i}\frac{\partial }{\partial
x_{n+i}}+X^{2n+i}\frac{\partial }{\partial x_{2n+i}}+X^{3n+i}\frac{\partial 
}{\partial x_{3n+i}},\,  \label{3.1}
\end{equation}%
where $X^{i}=\overset{.}{x_{i}},X^{n+i}=\overset{.}{x}_{n+i},X^{2n+i}=%
\overset{.}{x}_{2n+i},X^{3n+i}=\overset{.}{x}_{3n+i}$ and the dot indicates
the derivative with respect to time $t$. The vector field defined by%
\begin{equation*}
V_{F}(L)=F(X)=aX^{i}e^{\lambda }\frac{\partial L}{\partial x_{n+i}}%
-aX^{n+i}e^{-\lambda }\frac{\partial L}{\partial x_{i}}+aX^{2n+i}e^{\lambda }%
\frac{\partial L}{\partial x_{3n+i}}-aX^{3n+i}e^{-\lambda }\frac{\partial L}{%
\partial x_{2n+i}}
\end{equation*}%
is named a conformal \textit{Liouville vector field} on the
generalized-quaternionic K\"{a}hler manifold $(M,g,V)$. For $F$ , the closed
generalized-quaternionic K\"{a}hler form is the closed 2-form given by $\Phi
_{L}^{F}=-dd_{_{F}}L$ such that%
\begin{equation}
d_{_{F}}L=ae^{\lambda }\frac{\partial L}{\partial x_{n+i}}%
dx_{i}-ae^{-\lambda }\frac{\partial L}{\partial x_{i}}dx_{n+i}+ae^{\lambda }%
\frac{\partial L}{\partial x_{3n+i}}dx_{2n+i}-ae^{-\lambda }\frac{\partial L%
}{\partial x_{2n+i}}d_{3n+i}:\mathcal{F}(M)\rightarrow \wedge ^{1}{}M.
\label{3.11}
\end{equation}%
Then we have 
\begin{equation*}
\begin{array}{c}
\Phi _{L}^{F}=ae^{\lambda }\frac{\partial \lambda }{\partial x_{j}}\frac{%
\partial L}{\partial x_{n+i}}dx_{i}\wedge dx_{j}+ae^{\lambda }\frac{\partial
^{2}L}{\partial x_{j}\partial x_{n+i}}dx_{i}\wedge dx_{j}+ae^{-\lambda }%
\frac{\partial \lambda }{\partial x_{j}}\frac{\partial L}{\partial x_{i}}%
dx_{n+i}\wedge dx_{j} \\ 
-ae^{-\lambda }\frac{\partial ^{2}L}{\partial x_{j}\partial x_{i}}%
dx_{n+i}\wedge dx_{j}+ae^{\lambda }\frac{\partial \lambda }{\partial x_{j}}%
\frac{\partial L}{\partial x_{3n+i}}dx_{2n+i}\wedge dx_{j}+ae^{\lambda }%
\frac{\partial ^{2}L}{\partial x_{j}\partial x_{3n+i}}dx_{2n+i}\wedge dx_{j}
\\ 
+ae^{-\lambda }\frac{\partial \lambda }{\partial x_{j}}\frac{\partial L}{%
\partial x_{2n+i}}dx_{3n+i}\wedge dx_{j}-ae^{-\lambda }\frac{\partial ^{2}L}{%
\partial x_{j}\partial x_{2n+i}}dx_{3n+i}\wedge dx_{j}+ae^{\lambda }\frac{%
\partial \lambda }{\partial x_{n+j}}\frac{\partial L}{\partial x_{n+i}}%
dx_{i}\wedge dx_{n+j} \\ 
+ae^{\lambda }\frac{\partial ^{2}L}{\partial x_{n+j}\partial x_{n+i}}%
dx_{i}\wedge dx_{n+j}+ae^{-\lambda }\frac{\partial \lambda }{\partial x_{n+j}%
}\frac{\partial L}{\partial x_{i}}dx_{n+i}\wedge dx_{n+j}-ae^{-\lambda }%
\frac{\partial ^{2}L}{\partial x_{n+j}\partial x_{i}}dx_{n+i}\wedge dx_{n+j}
\\ 
+ae^{\lambda }\frac{\partial \lambda }{\partial x_{n+j}}\frac{\partial L}{%
\partial x_{3n+i}}dx_{2n+i}\wedge dx_{n+j}+ae^{\lambda }\frac{\partial ^{2}L%
}{\partial x_{n+j}\partial x_{3n+i}}dx_{2n+i}\wedge dx_{n+j}+ae^{-\lambda }%
\frac{\partial \lambda }{\partial x_{n+j}}\frac{\partial L}{\partial x_{2n+i}%
}dx_{3n+i}\wedge dx_{n+j} \\ 
-ae^{-\lambda }\frac{\partial ^{2}L}{\partial x_{n+j}\partial x_{2n+i}}%
dx_{3n+i}\wedge dx_{n+j}+ae^{\lambda }\frac{\partial \lambda }{\partial
x_{2n+j}}\frac{\partial L}{\partial x_{n+i}}dx_{i}\wedge
dx_{2n+j}+ae^{\lambda }\frac{\partial ^{2}L}{\partial x_{2n+j}\partial
x_{n+i}}dx_{i}\wedge dx_{2n+j}%
\end{array}%
\end{equation*}%
\begin{equation}
\begin{array}{c}
+ae^{-\lambda }\frac{\partial \lambda }{\partial x_{2n+j}}\frac{\partial L}{%
\partial x_{i}}dx_{n+i}\wedge dx_{2n+j}-ae^{-\lambda }\frac{\partial ^{2}L}{%
\partial x_{2n+j}\partial x_{i}}dx_{n+i}\wedge dx_{2n+j}+ae^{\lambda }\frac{%
\partial \lambda }{\partial x_{2n+j}}\frac{\partial L}{\partial x_{3n+i}}%
dx_{2n+i}\wedge dx_{2n+j} \\ 
+ae^{\lambda }\frac{\partial ^{2}L}{\partial x_{2n+j}\partial x_{3n+i}}%
dx_{2n+i}\wedge dx_{2n+j}+ae^{-\lambda }\frac{\partial \lambda }{\partial
x_{2n+j}}\frac{\partial L}{\partial x_{2n+i}}dx_{3n+i}\wedge
dx_{2n+j}-ae^{-\lambda }\frac{\partial ^{2}L}{\partial x_{2n+j}\partial
x_{2n+i}}dx_{3n+i}\wedge dx_{2n+j} \\ 
+ae^{\lambda }\frac{\partial \lambda }{\partial x_{3n+j}}\frac{\partial L}{%
\partial x_{n+i}}dx_{i}\wedge dx_{3n+j}+ae^{\lambda }\frac{\partial ^{2}L}{%
\partial x_{3n+j}\partial x_{n+i}}dx_{i}\wedge dx_{3n+j}+ae^{-\lambda }\frac{%
\partial \lambda }{\partial x_{3n+j}}\frac{\partial L}{\partial x_{i}}%
dx_{n+i}\wedge dx_{3n+j} \\ 
-ae^{-\lambda }\frac{\partial ^{2}L}{\partial x_{3n+j}\partial x_{i}}%
dx_{n+i}\wedge dx_{3n+j}+ae^{\lambda }\frac{\partial \lambda }{\partial
x_{3n+j}}\frac{\partial L}{\partial x_{3n+i}}dx_{2n+i}\wedge
dx_{3n+j}+ae^{\lambda }\frac{\partial ^{2}L}{\partial x_{3n+j}\partial
x_{3n+i}}dx_{2n+i}\wedge dx_{3n+j} \\ 
+ae^{-\lambda }\frac{\partial \lambda }{\partial x_{3n+j}}\frac{\partial L}{%
\partial x_{2n+i}}dx_{3n+i}\wedge dx_{3n+j}-ae^{-\lambda }\frac{\partial
^{2}L}{\partial x_{3n+j}\partial x_{2n+i}}dx_{3n+i}\wedge dx_{3n+j}%
\end{array}
\label{3.12}
\end{equation}%
Then we calculate%
\begin{equation*}
\begin{array}{c}
i_{X}\Phi _{L}^{F}=aX^{i}e^{\lambda }\frac{\partial \lambda }{\partial x_{j}}%
\frac{\partial L}{\partial x_{n+i}}dx_{j}-aX^{i}e^{\lambda }\frac{\partial
\lambda }{\partial x_{j}}\frac{\partial L}{\partial x_{n+i}}%
dx_{i}+aX^{i}e^{\lambda }\frac{\partial ^{2}L}{\partial x_{j}\partial x_{n+i}%
}dx_{j} \\ 
-aX^{i}e^{\lambda }\frac{\partial ^{2}L}{\partial x_{j}\partial x_{n+i}}%
dx_{i}-aX^{i}e^{-\lambda }\frac{\partial \lambda }{\partial x_{j}}\frac{%
\partial L}{\partial x_{i}}dx_{n+i}+aX^{i}e^{-\lambda }\frac{\partial ^{2}L}{%
\partial x_{j}\partial x_{i}}dx_{n+i} \\ 
-aX^{i}e^{\lambda }\frac{\partial \lambda }{\partial x_{j}}\frac{\partial L}{%
\partial x_{3n+i}}dx_{2n+i}-aX^{i}e^{\lambda }\frac{\partial ^{2}L}{\partial
x_{j}\partial x_{3n+i}}dx_{2n+i}-aX^{i}e^{-\lambda }\frac{\partial \lambda }{%
\partial x_{j}}\frac{\partial L}{\partial x_{2n+i}}dx_{3n+i} \\ 
+aX^{i}e^{-\lambda }\frac{\partial ^{2}L}{\partial x_{j}\partial x_{2n+i}}%
dx_{3n+i}+aX^{i}e^{\lambda }\frac{\partial \lambda }{\partial x_{n+j}}\frac{%
\partial L}{\partial x_{n+i}}dx_{n+j}+aX^{i}e^{\lambda }\frac{\partial ^{2}L%
}{\partial x_{n+j}\partial x_{n+i}}dx_{n+j} \\ 
+aX^{i}e^{\lambda }\frac{\partial \lambda }{\partial x_{2n+j}}\frac{\partial
L}{\partial x_{n+i}}dx_{2n+j}+aX^{i}e^{\lambda }\frac{\partial ^{2}L}{%
\partial x_{2n+j}\partial x_{n+i}}dx_{2n+j}+aX^{i}e^{\lambda }\frac{\partial
\lambda }{\partial x_{3n+j}}\frac{\partial L}{\partial x_{n+i}}dx_{3n+j} \\ 
+aX^{i}e^{\lambda }\frac{\partial ^{2}L}{\partial x_{3n+j}\partial x_{n+i}}%
dx_{3n+j}+aX^{n+i}e^{-\lambda }\frac{\partial \lambda }{\partial x_{j}}\frac{%
\partial L}{\partial x_{i}}dx_{j}-aX^{n+i}e^{-\lambda }\frac{\partial ^{2}L}{%
\partial x_{j}\partial x_{i}}dx_{j} \\ 
-aX^{n+i}e^{\lambda }\frac{\partial \lambda }{\partial x_{n+j}}\frac{%
\partial L}{\partial x_{n+i}}dx_{i}-aX^{n+i}e^{\lambda }\frac{\partial ^{2}L%
}{\partial x_{n+j}\partial x_{n+i}}dx_{i}+aX^{n+i}e^{-\lambda }\frac{%
\partial \lambda }{\partial x_{n+j}}\frac{\partial L}{\partial x_{i}}dx_{n+j}
\\ 
-aX^{n+i}e^{-\lambda }\frac{\partial \lambda }{\partial x_{n+j}}\frac{%
\partial L}{\partial x_{i}}dx_{n+i}-aX^{n+i}e^{-\lambda }\frac{\partial ^{2}L%
}{\partial x_{n+j}\partial x_{i}}dx_{n+j}+aX^{n+i}e^{-\lambda }\frac{%
\partial ^{2}L}{\partial x_{n+j}\partial x_{i}}dx_{n+i} \\ 
-aX^{n+i}e^{\lambda }\frac{\partial \lambda }{\partial x_{n+j}}\frac{%
\partial L}{\partial x_{3n+i}}dx_{2n+i}-aX^{n+i}e^{\lambda }\frac{\partial
^{2}L}{\partial x_{n+j}\partial x_{3n+i}}dx_{2n+i}-aX^{n+i}e^{-\lambda }%
\frac{\partial \lambda }{\partial x_{n+j}}\frac{\partial L}{\partial x_{2n+i}%
}dx_{3n+i} \\ 
+aX^{n+i}e^{-\lambda }\frac{\partial ^{2}L}{\partial x_{n+j}\partial x_{2n+i}%
}dx_{3n+i}+aX^{n+i}e^{-\lambda }\frac{\partial \lambda }{\partial x_{2n+j}}%
\frac{\partial L}{\partial x_{i}}dx_{2n+j}-aX^{n+i}e^{-\lambda }\frac{%
\partial ^{2}L}{\partial x_{2n+j}\partial x_{i}}dx_{2n+j} \\ 
+aX^{n+i}e^{-\lambda }\frac{\partial \lambda }{\partial x_{3n+j}}\frac{%
\partial L}{\partial x_{i}}dx_{3n+j}-aX^{n+i}e^{-\lambda }\frac{\partial
^{2}L}{\partial x_{3n+j}\partial x_{i}}dx_{3n+j}+aX^{2n+i}e^{\lambda }\frac{%
\partial \lambda }{\partial x_{j}}\frac{\partial L}{\partial x_{3n+i}}dx_{j}
\\ 
+aX^{2n+i}e^{\lambda }\frac{\partial ^{2}L}{\partial x_{j}\partial x_{3n+i}}%
dx_{j}+aX^{2n+i}e^{\lambda }\frac{\partial \lambda }{\partial x_{n+j}}\frac{%
\partial L}{\partial x_{3n+i}}dx_{n+j}+aX^{2n+i}e^{\lambda }\frac{\partial
^{2}L}{\partial x_{n+j}\partial x_{3n+i}}dx_{n+j}%
\end{array}%
\end{equation*}%
\begin{equation}
\begin{array}{c}
-aX^{2n+i}e^{\lambda }\frac{\partial \lambda }{\partial x_{2n+j}}\frac{%
\partial L}{\partial x_{n+i}}dx_{i}-aX^{2n+i}e^{\lambda }\frac{\partial ^{2}L%
}{\partial x_{2n+j}\partial x_{n+i}}dx_{i}-aX^{2n+i}e^{-\lambda }\frac{%
\partial \lambda }{\partial x_{2n+j}}\frac{\partial L}{\partial x_{i}}%
dx_{n+i} \\ 
+aX^{2n+i}e^{-\lambda }\frac{\partial ^{2}L}{\partial x_{2n+j}\partial x_{i}}%
dx_{n+i}+aX^{2n+i}e^{\lambda }\frac{\partial \lambda }{\partial x_{2n+j}}%
\frac{\partial L}{\partial x_{3n+i}}dx_{2n+j}-aX^{2n+i}e^{\lambda }\frac{%
\partial \lambda }{\partial x_{2n+j}}\frac{\partial L}{\partial x_{3n+i}}%
dx_{2n+i} \\ 
+aX^{2n+i}e^{\lambda }\frac{\partial ^{2}L}{\partial x_{2n+j}\partial
x_{3n+i}}dx_{2n+j}-aX^{2n+i}e^{\lambda }\frac{\partial ^{2}L}{\partial
x_{2n+j}\partial x_{3n+i}}dx_{2n+i}-aX^{2n+i}e^{-\lambda }\frac{\partial
\lambda }{\partial x_{2n+j}}\frac{\partial L}{\partial x_{2n+i}}dx_{3n+i} \\ 
+aX^{2n+i}e^{-\lambda }\frac{\partial ^{2}L}{\partial x_{2n+j}\partial x_{i}}%
dx_{3n+i}+aX^{2n+i}e^{\lambda }\frac{\partial \lambda }{\partial x_{3n+j}}%
\frac{\partial L}{\partial x_{3n+i}}dx_{3n+j}+aX^{2n+i}e^{\lambda }\frac{%
\partial ^{2}L}{\partial x_{3n+j}\partial x_{3n+i}}dx_{3n+j} \\ 
+aX^{3n+i}e^{-\lambda }\frac{\partial \lambda }{\partial x_{j}}\frac{%
\partial L}{\partial x_{2n+i}}dx_{j}-aX^{3n+i}e^{-\lambda }\frac{\partial
^{2}L}{\partial x_{j}\partial x_{2n+i}}dx_{j}+aX^{3n+i}e^{-\lambda }\frac{%
\partial \lambda }{\partial x_{n+j}}\frac{\partial L}{\partial x_{2n+i}}%
dx_{n+j} \\ 
-aX^{3n+i}e^{-\lambda }\frac{\partial ^{2}L}{\partial x_{n+j}\partial
x_{2n+i}}dx_{n+j}+aX^{3n+i}e^{-\lambda }\frac{\partial \lambda }{\partial
x_{2n+j}}\frac{\partial L}{\partial x_{2n+i}}dx_{2n+j}-aX^{3n+i}e^{-\lambda }%
\frac{\partial ^{2}L}{\partial x_{2n+j}\partial x_{i}}dx_{2n+j} \\ 
-aX^{3n+i}e^{\lambda }\frac{\partial \lambda }{\partial x_{3n+j}}\frac{%
\partial L}{\partial x_{n+i}}dx_{i}-aX^{3n+i}e^{\lambda }\frac{\partial ^{2}L%
}{\partial x_{3n+j}\partial x_{n+i}}dx_{i}-aX^{3n+i}e^{-\lambda }\frac{%
\partial \lambda }{\partial x_{3n+j}}\frac{\partial L}{\partial x_{i}}%
dx_{n+i} \\ 
+aX^{3n+i}e^{-\lambda }\frac{\partial ^{2}L}{\partial x_{3n+j}\partial x_{i}}%
dx_{n+i}-aX^{3n+i}e^{\lambda }\frac{\partial \lambda }{\partial x_{3n+j}}%
\frac{\partial L}{\partial x_{3n+i}}dx_{2n+i}-aX^{3n+i}e^{\lambda }\frac{%
\partial ^{2}L}{\partial x_{3n+j}\partial x_{3n+i}}dx_{2n+i} \\ 
+aX^{3n+i}e^{-\lambda }\frac{\partial \lambda }{\partial x_{3n+j}}\frac{%
\partial L}{\partial x_{2n+i}}dx_{3n+j}-aX^{3n+i}e^{-\lambda }\frac{\partial
\lambda }{\partial x_{3n+j}}\frac{\partial L}{\partial x_{2n+i}}dx_{3n+i} \\ 
-aX^{3n+i}e^{-\lambda }\frac{\partial ^{2}L}{\partial x_{3n+j}\partial
x_{2n+i}}dx_{3n+j}+aX^{3n+i}e^{-\lambda }dx_{3n+i}%
\end{array}
\label{3.13}
\end{equation}%
Energy function is%
\begin{equation}
\begin{array}{c}
E_{L}^{F}=V_{F}(L)-L=aX^{i}e^{\lambda }\frac{\partial L}{\partial x_{n+i}}%
-aX^{n+i}e^{-\lambda }\frac{\partial L}{\partial x_{i}}+aX^{2n+i}e^{\lambda }%
\frac{\partial L}{\partial x_{3n+i}}-aX^{3n+i}e^{-\lambda }\frac{\partial L}{%
\partial x_{2n+i}}-L%
\end{array}
\label{3.14}
\end{equation}%
and hence%
\begin{equation*}
\begin{array}{c}
dE_{L}^{F}=aX^{i}e^{\lambda }\frac{\partial \lambda }{\partial x_{j}}\frac{%
\partial L}{\partial x_{n+i}}dx_{j}+aX^{i}e^{\lambda }\frac{\partial ^{2}L}{%
\partial x_{j}\partial x_{n+i}}dx_{j}+aX^{n+i}e^{-\lambda }\frac{\partial
\lambda }{\partial x_{j}}\frac{\partial L}{\partial x_{i}}dx_{j} \\ 
-aX^{n+i}e^{-\lambda }\frac{\partial ^{2}L}{\partial x_{j}\partial x_{i}}%
dx_{j}+aX^{2n+i}e^{\lambda }\frac{\partial \lambda }{\partial x_{j}}\frac{%
\partial L}{\partial x_{3n+i}}dx_{j}+aX^{2n+i}e^{\lambda }\frac{\partial
^{2}L}{\partial x_{j}\partial x_{3n+i}}dx_{j} \\ 
+aX^{3n+i}e^{-\lambda }\frac{\partial \lambda }{\partial x_{j}}\frac{%
\partial L}{\partial x_{2n+i}}dx_{j}-aX^{3n+i}e^{-\lambda }\frac{\partial
^{2}L}{\partial x_{j}\partial x_{2n+i}}dx_{j}-\frac{\partial L}{\partial
x_{j}}dx_{j} \\ 
+aX^{i}e^{\lambda }\frac{\partial \lambda }{\partial x_{n+j}}\frac{\partial L%
}{\partial x_{n+i}}dx_{n+j}+aX^{i}e^{\lambda }\frac{\partial ^{2}L}{\partial
x_{n+j}\partial x_{n+i}}dx_{n+j}+aX^{n+i}e^{-\lambda }\frac{\partial \lambda 
}{\partial x_{n+j}}\frac{\partial L}{\partial x_{i}}dx_{n+j} \\ 
-aX^{n+i}e^{-\lambda }\frac{\partial ^{2}L}{\partial x_{n+j}\partial x_{i}}%
dx_{n+j}+aX^{2n+i}e^{\lambda }\frac{\partial \lambda }{\partial x_{n+j}}%
\frac{\partial L}{\partial x_{3n+i}}dx_{n+j}+aX^{2n+i}e^{\lambda }\frac{%
\partial ^{2}L}{\partial x_{n+j}\partial x_{3n+i}}dx_{n+j} \\ 
+aX^{3n+i}e^{-\lambda }\frac{\partial \lambda }{\partial x_{n+j}}\frac{%
\partial L}{\partial x_{2n+i}}dx_{n+j}-aX^{3n+i}e^{-\lambda }\frac{\partial
^{2}L}{\partial x_{n+j}\partial x_{2n+i}}dx_{n+j}-\frac{\partial L}{\partial
x_{n+j}}dx_{n+j}%
\end{array}%
\end{equation*}%
\begin{equation}
\begin{array}{c}
+aX^{i}e^{\lambda }\frac{\partial \lambda }{\partial x_{2n+j}}\frac{\partial
L}{\partial x_{n+i}}dx_{2n+j}+aX^{i}e^{\lambda }\frac{\partial ^{2}L}{%
\partial x_{2n+j}\partial x_{n+i}}dx_{2n+j}+aX^{n+i}e^{-\lambda }\frac{%
\partial \lambda }{\partial x_{2n+j}}\frac{\partial L}{\partial x_{i}}%
dx_{2n+j} \\ 
-aX^{n+i}e^{-\lambda }\frac{\partial ^{2}L}{\partial x_{2n+j}\partial x_{i}}%
dx_{2n+j}+aX^{2n+i}e^{\lambda }\frac{\partial \lambda }{\partial x_{2n+j}}%
\frac{\partial L}{\partial x_{3n+i}}dx_{2n+j}+aX^{2n+i}e^{\lambda }\frac{%
\partial ^{2}L}{\partial x_{2n+j}\partial x_{3n+i}}dx_{2n+j} \\ 
+aX^{3n+i}e^{-\lambda }\frac{\partial \lambda }{\partial x_{2n+j}}\frac{%
\partial L}{\partial x_{2n+i}}dx_{2n+j}-aX^{3n+i}e^{-\lambda }\frac{\partial
^{2}L}{\partial x_{2n+j}\partial x_{2n+i}}dx_{2n+j}-\frac{\partial L}{%
\partial x_{2n+j}}dx_{2n+j} \\ 
+aX^{i}e^{\lambda }\frac{\partial \lambda }{\partial x_{3n+j}}\frac{\partial
L}{\partial x_{n+i}}dx_{3n+j}+aX^{i}e^{\lambda }\frac{\partial ^{2}L}{%
\partial x_{3n+j}\partial x_{n+i}}dx_{3n+j}+aX^{n+i}e^{-\lambda }\frac{%
\partial \lambda }{\partial x_{3n+j}}\frac{\partial L}{\partial x_{i}}%
dx_{3n+j} \\ 
-aX^{n+i}e^{-\lambda }\frac{\partial ^{2}L}{\partial x_{3n+j}\partial x_{i}}%
dx_{3n+j}+aX^{2n+i}e^{\lambda }\frac{\partial \lambda }{\partial x_{3n+j}}%
\frac{\partial L}{\partial x_{3n+i}}dx_{3n+j}+aX^{2n+i}e^{\lambda }\frac{%
\partial ^{2}L}{\partial x_{3n+j}\partial x_{3n+i}}dx_{3n+j} \\ 
+aX^{3n+i}e^{-\lambda }\frac{\partial \lambda }{\partial x_{3n+j}}\frac{%
\partial L}{\partial x_{2n+i}}dx_{3n+j}-aX^{3n+i}e^{-\lambda }\frac{\partial
^{2}L}{\partial x_{3n+j}\partial x_{2n+i}}dx_{3n+j}-\frac{\partial L}{%
\partial x_{3n+j}}dx_{3n+j}%
\end{array}
\label{3.114}
\end{equation}%
Using (\ref{1.1}) and also considering an integral curve of $X,$ then we
obtain the equation given by%
\begin{equation*}
\begin{array}{c}
-aX^{i}e^{\lambda }\frac{\partial \lambda }{\partial x_{j}}\frac{\partial L}{%
\partial x_{n+i}}dx_{i}-aX^{i}e^{\lambda }\frac{\partial ^{2}L}{\partial
x_{j}\partial x_{n+i}}dx_{i}-aX^{n+i}e^{\lambda }\frac{\partial \lambda }{%
\partial x_{n+j}}\frac{\partial L}{\partial x_{n+i}}dx_{i} \\ 
-aX^{n+i}e^{\lambda }\frac{\partial ^{2}L}{\partial x_{n+j}\partial x_{n+i}}%
dx_{i}-aX^{2n+i}e^{\lambda }\frac{\partial \lambda }{\partial x_{2n+j}}\frac{%
\partial L}{\partial x_{n+i}}dx_{i}-aX^{2n+i}e^{\lambda }\frac{\partial ^{2}L%
}{\partial x_{2n+j}\partial x_{n+i}}dx_{i} \\ 
-aX^{3n+i}e^{\lambda }\frac{\partial \lambda }{\partial x_{3n+j}}\frac{%
\partial L}{\partial x_{n+i}}dx_{i}-aX^{3n+i}e^{\lambda }\frac{\partial ^{2}L%
}{\partial x_{3n+j}\partial x_{n+i}}dx_{i}+\frac{\partial L}{\partial x_{j}}%
dx_{j} \\ 
-aX^{i}e^{-\lambda }\frac{\partial \lambda }{\partial x_{j}}\frac{\partial L%
}{\partial x_{i}}dx_{n+i}+aX^{i}e^{-\lambda }\frac{\partial ^{2}L}{\partial
x_{j}\partial x_{i}}dx_{n+i}-aX^{n+i}e^{-\lambda }\frac{\partial ^{2}L}{%
\partial x_{n+j}\partial x_{i}}dx_{n+j} \\ 
+aX^{n+i}e^{-\lambda }\frac{\partial ^{2}L}{\partial x_{n+j}\partial x_{i}}%
dx_{n+i}-aX^{2n+i}e^{-\lambda }\frac{\partial \lambda }{\partial x_{2n+j}}%
\frac{\partial L}{\partial x_{i}}dx_{n+i}+aX^{2n+i}e^{-\lambda }\frac{%
\partial ^{2}L}{\partial x_{2n+j}\partial x_{i}}dx_{n+i} \\ 
-aX^{3n+i}e^{-\lambda }\frac{\partial \lambda }{\partial x_{3n+j}}\frac{%
\partial L}{\partial x_{i}}dx_{n+i}+aX^{3n+i}e^{-\lambda }\frac{\partial
^{2}L}{\partial x_{3n+j}\partial x_{i}}dx_{n+i}+\frac{\partial L}{\partial
x_{n+j}}dx_{n+j}%
\end{array}%
\end{equation*}%
\begin{equation}
\begin{array}{c}
-aX^{i}e^{\lambda }\frac{\partial \lambda }{\partial x_{j}}\frac{\partial L}{%
\partial x_{3n+i}}dx_{2n+i}-aX^{i}e^{\lambda }\frac{\partial ^{2}L}{\partial
x_{j}\partial x_{3n+i}}dx_{2n+i}-aX^{n+i}e^{\lambda }\frac{\partial \lambda 
}{\partial x_{n+j}}\frac{\partial L}{\partial x_{3n+i}}dx_{2n+i} \\ 
-aX^{n+i}e^{\lambda }\frac{\partial ^{2}L}{\partial x_{n+j}\partial x_{3n+i}}%
dx_{2n+i}+aX^{2n+i}e^{\lambda }\frac{\partial ^{2}L}{\partial
x_{2n+j}\partial x_{3n+i}}dx_{2n+j}-aX^{2n+i}e^{\lambda }\frac{\partial ^{2}L%
}{\partial x_{2n+j}\partial x_{3n+i}}dx_{2n+i} \\ 
-aX^{3n+i}e^{\lambda }\frac{\partial \lambda }{\partial x_{3n+j}}\frac{%
\partial L}{\partial x_{3n+i}}dx_{2n+i}-aX^{3n+i}e^{\lambda }\frac{\partial
^{2}L}{\partial x_{3n+j}\partial x_{3n+i}}dx_{2n+i}+\frac{\partial L}{%
\partial x_{2n+j}}dx_{2n+j} \\ 
-aX^{i}e^{-\lambda }\frac{\partial \lambda }{\partial x_{j}}\frac{\partial L%
}{\partial x_{2n+i}}dx_{3n+i}+aX^{i}e^{-\lambda }\frac{\partial ^{2}L}{%
\partial x_{j}\partial x_{2n+i}}dx_{3n+i}-aX^{n+i}e^{-\lambda }\frac{%
\partial \lambda }{\partial x_{n+j}}\frac{\partial L}{\partial x_{2n+i}}%
dx_{3n+i} \\ 
+aX^{n+i}e^{-\lambda }\frac{\partial ^{2}L}{\partial x_{n+j}\partial x_{2n+i}%
}dx_{3n+i}-aX^{2n+i}e^{-\lambda }\frac{\partial \lambda }{\partial x_{2n+j}}%
\frac{\partial L}{\partial x_{2n+i}}dx_{3n+i}+aX^{2n+i}e^{-\lambda }\frac{%
\partial ^{2}L}{\partial x_{2n+j}\partial x_{i}}dx_{3n+i} \\ 
-aX^{3n+i}e^{-\lambda }\frac{\partial ^{2}L}{\partial x_{3n+j}\partial
x_{2n+i}}dx_{3n+j}+aX^{3n+i}e^{-\lambda }dx_{3n+i}+\frac{\partial L}{%
\partial x_{3n+j}}dx_{3n+j}=0%
\end{array}
\label{3.15}
\end{equation}%
Then we have the equations%
\begin{equation}
\begin{array}{l}
a\frac{\partial }{\partial t}\left( e^{-\lambda }\frac{\partial L}{\partial
x_{i}}\right) +\frac{\partial L}{\partial x_{n+i}}=0,\text{ }a\frac{\partial 
}{\partial t}\left( e^{\lambda }\frac{\partial L}{\partial x_{n+i}}\right) -%
\frac{\partial L}{\partial x_{i}}=0, \\ 
a\,\frac{\partial }{\partial t}\left( e^{-\lambda }\frac{\partial L}{%
\partial x_{2n+i}}\right) +\frac{\partial L}{\partial x_{3n+i}}=0,\text{\ }a%
\frac{\partial }{\partial t}\left( e^{\lambda }\frac{\partial L}{\partial
x_{3n+i}}\right) -\frac{\partial L}{\partial x_{2n+i}}=0,%
\end{array}
\label{3.2}
\end{equation}%
such that the equations calculated in (\ref{3.2}) are named \textit{%
Euler-Lagrange equations} constructed on a generalized-quaternionic K\"{a}%
hler manifold $(M,g,V)$ by means of $\Phi _{L}^{F}$ and thus the triple $%
(M,\Phi _{L}^{F},X)$ is called a \textit{mechanical system }on a
generalized-quaternionic K\"{a}hler manifold $(M,g,V)$\textit{. }Considering
the above operations and also taking the following vector fields 
\begin{equation}
Y=Y^{i}\frac{\partial }{\partial x_{i}}+Y^{n+i}\frac{\partial }{\partial
x_{n+i}}+Y^{2n+i}\frac{\partial }{\partial x_{2n+i}}+Y^{3n+i}\frac{\partial 
}{\partial x_{3n+i}},  \label{3.3}
\end{equation}%
\begin{equation}
Z=Z^{i}\frac{\partial }{\partial x_{i}}+Z^{n+i}\frac{\partial }{\partial
x_{n+i}}+Z^{2n+i}\frac{\partial }{\partial x_{2n+i}}+Z^{3n+i}\frac{\partial 
}{\partial x_{3n+i}},  \label{3.4}
\end{equation}%
we obtain the following Euler-Lagrange equations for quantum and classical
mechanics by means of $\Phi _{L}^{G}$ and $\Phi _{L}^{H}$ on a
generalized-quaternionic K\"{a}hler manifold $(M,g,V),$ respectively:%
\begin{equation}
\begin{array}{l}
b\frac{\partial }{\partial t}\left( e^{-\lambda }\frac{\partial L}{\partial
x_{i}}\right) +\frac{\partial L}{\partial x_{2n+i}}=0,\text{ \ }b\frac{%
\partial }{\partial t}\left( e^{-\lambda }\frac{\partial L}{\partial x_{n+i}}%
\right) -\frac{\partial L}{\partial x_{3n+i}}=0, \\ 
b\frac{\partial }{\partial t}\left( e^{\lambda }\frac{\partial L}{\partial
x_{2n+i}}\right) +\frac{\partial L}{\partial x_{i}}=0,\text{ \ }b\frac{%
\partial }{\partial t}\left( e^{\lambda }\frac{\partial L}{\partial x_{3n+i}}%
\right) -\frac{\partial L}{\partial x_{n+i}}=0.%
\end{array}
\label{3.5}
\end{equation}%
\begin{equation}
\begin{array}{l}
ab\frac{\partial }{\partial t}\left( e^{-\lambda }\frac{\partial L}{\partial
x_{i}}\right) +\frac{\partial L}{\partial x_{3n+i}}=0,\text{ \ }ab\frac{%
\partial }{\partial t}\left( e^{-\lambda }\frac{\partial L}{\partial x_{n+i}}%
\right) +\frac{\partial L}{\partial x_{2n+i}}=0, \\ 
ab\frac{\partial }{\partial t}\left( e^{\lambda }\frac{\partial L}{\partial
x_{2n+i}}\right) +\frac{\partial L}{\partial x_{n+i}}=0,\text{ \ }ab\frac{%
\partial }{\partial t}\left( e^{\lambda }\frac{\partial L}{\partial x_{3n+i}}%
\right) +\frac{\partial L}{\partial x_{i}}=0.%
\end{array}
\label{3.6}
\end{equation}%
Thus the equations introduced by (\ref{3.5}) and (\ref{3.6}) infer Conformal 
\textit{Euler-Lagrange equations} constructed by means of $\Phi _{L}^{G}$
and $\Phi _{L}^{H}$ on a generalized-quaternionic K\"{a}hler manifold $%
(M,g,V)$ and then the triples $(M,\Phi _{L}^{G},X)$ and $(M,\Phi _{L}^{H},X)$
are named \textit{mechanical systems }on a generalized-quaternionic K\"{a}%
hler manifold $(M,g,V)$. \textbf{Hence the equations found by (\ref{3.2},\ref%
{3.5},\ref{3.6}) easily seen extremely useful in applications from
Euler-Lagrangian Mechanics, Quantum Physics, Optimal Control , Biology and
Fluid Dynamics \cite{miron,weyl}.}

\section{Conclusion}

Given the above equations, Euler-Lagrangian mechanical systems have
intrinsically been defined taking into account a canonical local basis $%
\{F,G,H\}$ of $V$ \ that they are defined on a generalized-quaternionic K%
\"{a}hler manifold $(M,g,V).$ The paths of semisprays $X,G,H$ on the
generalized-quaternionic K\"{a}hler manifold are the solutions
Euler-Lagrange equations raised in (\ref{3.2}), (\ref{3.5}) and (\ref{3.6}).
These equations are introduced by a canonical local basis $\{F,G,H\}$ of
vector bundle $V$ on a generalized-quaternionic K\"{a}hler manifold $(M,g,V)$%
. If this equations $a=1$ and $b=1$ are selected, the equations given in (%
\ref{3.2}), (\ref{3.5}) and (\ref{3.6}) we say to be Euler-Lagrange
equations on a quaternionic K\"{a}hler manifold. If $a=1$ and $b=-1$, the
equations given in (\ref{3.2}), (\ref{3.5}) and (\ref{3.6}) we say to be
Euler-Lagrange equations on a para-quaternionic K\"{a}hler manifold. In
these days, Lagrangian models arise to be a very important tool since they
present a simple method to describe the model for mechanical systems. One
can be proved that the obtained equations are very important to explain the
rotational spatial mechanical-physical problems. For this reason, the found
equations are only considered to be a first step to realize how a
generalized-quaternionic geometry has been used in solving problems in
different physical area.

In the literature, the equations, which explains the linear orbits of the
objects, were presented. This study explained the non-linear orbits of the
objects in the space by the help of revised equations using conformal
structure.

Our proposal for future research, the Lagrange mechanical equations derived
on a generalized-quaternionic K\"{a}hler manifold are suggested to deal with
problems in electrical, magnetical and gravitational fields of quantum and
classical mechanics of physics.

\end{document}